\documentclass{article}

\usepackage[mathcal]{eucal} \usepackage{amsmath} \usepackage{amsfonts}

\usepackage{graphicx} \usepackage{url} \usepackage[square]{natbib}
\usepackage{subfigure}

\newcommand{\ensR}{\mathbb{R}} \newcommand{\Forall}{\forall \hspace*{0.7mm}}
\newcommand{\code}[1]{\texttt{#1}}

\begin{document}

\title{Modeling Wildland Fire Propagation\\ with Level Set Methods}

\author{V. Mallet$^1$, D. E. Keyes$^2$, F.  E. Fendell$^3$\\
  {\small$^1$ \'Ecole des Ponts, Marne-la-Vall\'ee, France}\\
  {\small$^2$ Columbia University, City of New York, USA}\\
  {\small$^3$ Northrop Grumman Space Technology, Redondo Beach, CA, USA}}

\date{August 2007}

\maketitle

{\it Level set methods are versatile and extensible techniques for general
  front tracking problems, including the practically important problem of
  predicting the advance of a firefront across expanses of surface vegetation.
  Given a rule, empirical or otherwise, to specify the rate of advance of an
  infinitesimal segment of firefront arc normal to itself (i.e., given the
  firespread rate as a function of known local parameters relating to
  topography, vegetation, and meteorology), level set methods harness the well
  developed mathematical machinery of hyperbolic conservation laws on Eulerian
  grids to evolve the position of the front in time.  Topological challenges
  associated with the swallowing of islands and the merger of fronts are
  tractable.

  The principal goals of this paper are to: collect key results from the two
  largely distinct scientific literatures of level sets and firespread;
  demonstrate the practical value of level set methods to wildland fire
  modeling through numerical experiments; probe and address current
  limitations; and propose future directions in the simulation of, and the
  development of decision-aiding tools to assess countermeasure options for,
  wildland fires.  In addition, we introduce a freely available
  two-dimensional level set code used to produce the numerical results of this
  paper and designed to be extensible to more complicated configurations. \\ }
  
Keywords: Wildland firespread, level set methods, Multivac software

\section{Introduction}

Wildland fire modeling has received attention for decades, due to the
sometimes disastrous consequences of large fires, and the tremendous costs of
often ineffectual, possibly even counterproductive firefighting
\citep{pyne04tending}.  For the practically important scenario of wind-aided
firespread, one seeks a computationally efficient model, useful not only
offline (for pre-crisis planning, e.g., placement of access roads, firebreaks,
and reservoirs, and scoping of fuel-reduction burning, and post-crisis review,
e.g., personnel training, litigation), but also during a crisis (i.e.,
real-time guidance for evacuation and firefighting).  For computational
efficiency, such that the benefits of ensemble forecasting
\citep{palmer05representing} are readily accessible from a model, advantage
should be taken of the inherent scale separation of: (1) the
kilometer-and-larger, landscape-dominated scales of the local atmospheric
dynamics; and (2) the one-meter-and-smaller scales of the local combustion
dynamics.  Even with advanced techniques and access to exceptional
contemporary computing facilities, numerical simulations (of turbulent flows)
that proceed from fundamental principles are challenged to resolve accurately
in real time phenomena with spatial scales spanning much more than two orders
of magnitude [H.R. Baum, private communication].  Thus, the feasibility of a
direct numerical simulation encompassing the multivaried processes of wildland
fire propagation \citep{coen03wildfire} may be decades off
\citep{jenkins01forest}.  Moreover, at least many attempts (albeit usually
problematic) at parameterization of subgridscale phenomena in terms of
gridscale variables have been undertaken by meteorologists for cumulus
convection, turbulent transport, and radiative transfer.  However,
meteorologists have extremely limited experience with the parameterization of
combustion dynamics for weather-dependent wildland firespread; even if such
parameterization be possible, it remains unknown.  Furthermore, data
collection in wildland fires is so piecemeal, irregular, and of uncertain
accuracy that, for many years to come, the data better suit reinitialization
of a simplistic model than assimilation into an ongoing calculation with a
highly detailed model.

Accordingly, in this study, attention is focused on a minimalist treatment of
the firefront, idealized as an interface between expanses of burned and
unburned vegetation.  This treatment is consistent with the typically limited,
only gross characterization available for the vegetation at issue, since the
vagaries of ignition events are difficult to anticipate, and maintaining an
updated inventory for the huge area of wildlands in (say) the USA is daunting.
This simplistic interfacial approach to the fire dynamics, easily executed in
minutes on a laptop given the requisite meteorological and other input fields,
reserves computational resources for the difficult, more critical, and mostly
yet-to-be-undertaken landscape-scale weather forecasting targeted for
real-time wildfire applications.

The upshot is that simple persistence models are adopted for the wind field
(and thermodynamic variables) in the study undertaken here.  Also, attention
is limited to a one-way interaction between the meteorology and the
firespread, though future extension to two-way interaction by use of an
iterative procedure may be envisioned.  Simplistic modeling still may provide
the key macroscopic fire behavior sufficiently accurately for practical
purposes (including estimates of smoke and pollutant generation), even for
circumstances for which the simplification is not formally justifiable.  In
fact, observational data of wind-aided firefront progression in wildland are
today typically sparse, so that not much more than the output of a simplistic
model can be meaningfully validated and tuned.  Moreover, the use of relevant
mathematical methods to perform model selection, to carry out efficient
parameter estimation, and to account for the uncertainty in predictions is
facilitated by focusing on less detailed models with fewer parameters.  In
this paper, we mainly address the first step, which is to achieve proper
forward simulations.

One of the most widely used models was devised by Rothermel
\citep{rothermel72mathematical} to predict the rate of firespread, with focus
on the head of a wind-aided fire.  Because predictions of the Rothermel
treatment have been found to be at odds with some observations, efforts to
improve this spatially one-dimensional semi-empirical treatment, and to
supplement the data upon which it is based, have been undertaken, especially
in recent years \citep{carlton03impact}.  Extension from a focus exclusively
on the head of the fire seeks to evolve the configuration of the entire fire
perimeter, possibly of multiple fire perimeters.  In this study, and
typically, the firefront, even a moderate fraction of an hour after a
localized ignition in fire-prone vegetation, is taken to be a closed curve
projected on a plane (the ground may not be flat).  Such simulations of
firespread have been performed \citep{finney98farsite} with the so-called
marker technique, which discretizes a front into a set of marker particles,
and advances the front through updates of the particle positions.
Parenthetically, as a problematic step, the updating by Finney takes each
marker on the front to evolve identically to an idealization of how a front
evolves from a single isolated ignition site in an unbounded expanse of
vegetation, in the presence of a wind.  In any case, even though applied
projects have supported software development \citep{finney98farsite}, still
from a computational point of view, only a few, largely equivalent
methodological developments have been undertaken
\citep[e.g.][]{andre06forest}.  In this paper, we apply {\it level set
  methods} \citep{osher88fronts,sethian99level} to calculate firefront
evolution.

In Section~\ref{sec:fire_model}, we introduce wildland firespread models,
especially a semi-empirical, equilibrium-type model proposed in
\cite{fendell01forest} for wind-aided firespread across surface-layer,
chaparral-type, burning-prone vegetation.  (In commonly adopted
equilibrium-type models, the firespread rate depends on only the parametric
values holding locally and instantaneously, so the firespread rate is taken to
adjust indefinitely rapidly to any temporal and spatial change.)
Section~\ref{sec:level_set} provides a brief introduction to level set
methods. Section~\ref{sec:code} describes the Multivac level set package that
has been applied in this paper to the firespread problem.  A quick description
of its performance is presented in Section~\ref{sec:complexity}.  Finally,
results of firespread simulations with different idealized environmental
conditions are reported in Section~\ref{sec:applying}.








\section{Front Propagation Functions for Wildland Fires}
\label{sec:fire_model}






Even if theory and/or measurement furnished complete, perfect knowledge of the
topography, vegetation, and meteorology at a site at a given time (e.g.,
furnished the locally pertinent values of all parameters in functional forms
capable of representing these three types of input), still one currently
possesses very incomplete, imperfect knowledge of the ``rules'' that would
yield the physically observed rate of firespread from the input.  Achieving
knowledge of firespread ``rules'' sufficiently accurate for practical purposes
may well lag emplacing means for observing and collecting exhaustive input
data.

As already noted, a fire-growth simulation such as FARSITE
\citep{finney98farsite} seems unlikely to reach its potential as long as it
seeks to describe the rate of firespread at all orientations to the direction
of the sustained low-level ambient wind from spread-rate modeling focused on
the direction of the wind \citep[e.g.][]{rothermel72mathematical}.  On the
other hand, posing a different rule for the spread rate at every possible
orientation to the wind defeats the goal of simplicity.

\subsection{Wind-aided wildland fire spread}

Fendell and Wolff \citep{fendell01forest} addressed this dilemma in developing
a model dedicated to wind-aided wildland fires that spread rapidly over level
terrain with dry, moderately sparse fuel, taken here to be uniformly
distributed to permit concentration on wind effects.  Parenthetically, for
consistency with modeling in which the firefront is idealized as an interface
moving according to a semi-empirical rule, only a minimal amount of
information about the surface-layer fuel is required, mainly the mass loading
consumed with firefront passage (``available''-fuel loading).

\begin{figure}
  \includegraphics[width=0.8\textwidth]{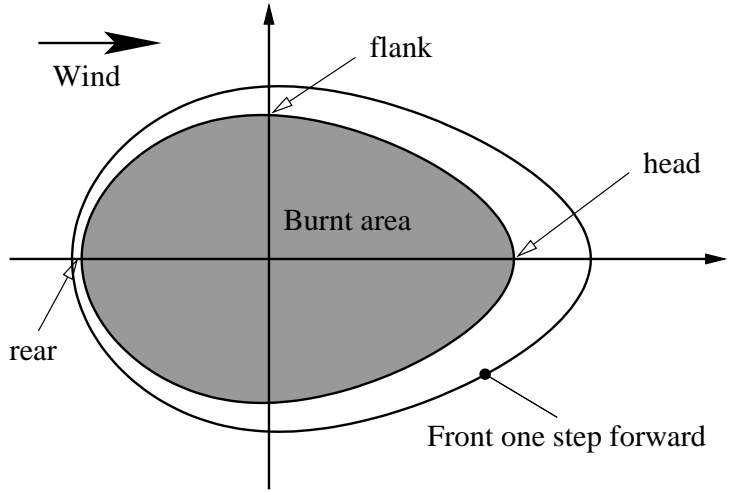}
  \caption{Fendell and Wolff model introduces velocities at the rear (against
    the wind), at the head (in the wind direction) and at the flanks.
    \citep{fendell01forest}}
  \label{fig:front}
\end{figure}

The Fendell and Wolff model focuses on front velocities at the rear of the
front (where propagation is against the wind), at the head of the front (where
propagation is with the wind), and on the flanks (where propagation is across
the wind direction) -- see Figure~\ref{fig:front}. The firespread velocities
primarily depend on the wind velocity $U$.  At the rear, the front advances
relatively slowly against the oncoming wind, since hot combustion products
tend to be blown over an already burned area.  The velocity at the rear is
denoted $\varepsilon(U)$.  At the head, the velocity $h(U)$ is relatively
large, since hot combustion products tend to be blown over a yet-to-burn area,
in which discrete fuel elements are heated toward ignition by
convective-conductive transfer.  Both analytic modeling and laboratory
experiments have shown that $h(U)$ is roughly proportional to $\sqrt{U}$
\citep{wolff91wind}.  At the flanks, the (spread-aiding) wind component along
the normal to the front is zero, but observationally the front advances faster
than in the absence of wind. As a speculation, a more meticulous treatment
would find that, at the nominal flank, the configuration is convoluted, and
firespread is alternately with and against the wind.  Of course, were the wind
direction constant, limiting attention to the head would seem adequate, but,
in fact, change in wind direction may (rapidly) result in an interchange of
the locations of the flank and head -- an interchange sometimes associated
with tragic consequence for firefighters.

The velocities (the terminology henceforth adopted, for brevity, in place of
firespread rates) proposed in \cite{fendell01forest} are
\begin{equation}
  \varepsilon(U) = \varepsilon_0 \exp(-\varepsilon_1 U), \quad f(U) =
  \varepsilon_0 + c_1 U \exp(-c_2 U), \quad h(U) = \varepsilon_0 + a \sqrt{U},
  \label{eq:velocities}
\end{equation}
where $\varepsilon_0$, $\varepsilon_1$, $c_1$, $c_2$ and $a$ are parameters
(with readily inferred dimensionality) depending on the mass loading of fuel
and other parameters characterizing the fuel bed, but independent of $U$.

The velocity is then provided at any point on the front through a
``trigonometric interpolation'':
\begin{equation}
  \begin{array}{l}
    F(U, \theta) = f(U \sin^m \theta) + h(U \cos^n \theta) \quad \textrm{if
      $|\theta| \leq \frac{\pi}{2}$},\\
    F(U, \theta) = f(U \sin^m \theta) + \varepsilon(U \cos^2 \theta) \quad
    \textrm{if $|\theta| > \frac{\pi}{2}$},
  \end{array}
  \label{eq:ext_velocities}
\end{equation}
where $\theta$ is the angle between the wind direction and the normal to the
front. We set $m = 2$. In this paper, parameter $n$ is set to $\frac{3}{2}$
and is significant since it determines the overall shape of the front from the
flanks to the head.

To summarize, the velocity is, for all $U\in\ensR_+$ and $ \theta \in ]-\pi,
\pi[$,
\begin{equation}
  \begin{array}{rl}
    F(U, \theta) = & \varepsilon_0 \sin^2 \theta + c_0 U \sin^2 \theta
    \hspace{0.5mm} \exp \left(-c_1 U \sin^2 \theta \right) \\ & + \left\{
      \begin{array}{llll}
        \displaystyle \varepsilon_0 \cos^2 \theta + a \sqrt{U}
        \cos^n \theta \quad & \textrm{if $|\theta| \leq
          \frac{\pi}{2}$}\\ \displaystyle \varepsilon_0 \cos^2 \theta
        \hspace{0.5mm} \exp \left(-\varepsilon_1 U \cos^2 \theta \right) \quad
        & \textrm{if $|\theta| > \frac{\pi}{2}$}
      \end{array}
    \right..
  \end{array}
  \label{eq:fire_model}
\end{equation}

\subsection{Simplified model}

Based on the numerical experiments carried out with the level set code
Multivac (Section~\ref{sec:code}), the model~(\ref{eq:fire_model}) proposed in
\cite{fendell01forest} has been modified. First, the parameter $n$ has been
set to $\frac{3}{2}$ instead of $1$. Second, the model has been simplified
without losing its main features, primarily the overall shape of the
firefront. The new model reads
\begin{equation}
  \begin{array}{ll}
    F(U, \theta) = \varepsilon_0 + c_1 \sqrt{U} \cos^n \theta & \quad
    \textrm{if $|\theta| \leq \frac{\pi}{2}$},\\
    F(U, \theta) = \varepsilon_0 (\alpha + (1 - \alpha)|\sin \theta|) & \quad
    \textrm{if $|\theta| > \frac{\pi}{2}$},
  \end{array}
  \label{eq:simplified_model}
\end{equation}
where $\alpha \in [0, 1]$ is the ratio between the velocity at the rear
($\alpha \varepsilon_0$) and the velocity at the flanks ($\varepsilon_0$).
Velocities at the rear and at the flanks no longer depend on the wind, since
their dependence on the wind speed is hard to model accurately and has little
impact on the overall front location. The velocity at the head is the same as
in the ``full'' model~(\ref{eq:fire_model}).

The simplified model is easier to tune, either via direct trials or with
systematic methods for parameter estimation (which may require derivatives of
the model with respect to its parameters). All results in this paper are for
the simplified model. However, results for the ``full'' model would appear
roughly the same.

\section{Level Set and Fast Marching Methods}
\label{sec:level_set}



First introduced in \citet{osher88fronts}, level set methods are Eulerian
schemes for tracking fronts propagating according to a given speed function.
In this section, we explain basic features of the level set methods used for
firespread modeling.

\subsection{Mathematical basis and technique}

\subsubsection{Definitions}

Assume the front evolves from the initial time $t = 0$ to the final time $t =
T_f$. For all $t \in \left[0, T_f\right]$, the front at time $t$ is the set of
points (in $\ensR^N$) $\Gamma(t)$. We define $\Gamma_0=\Gamma(0)$ as the
initial front.

For all $t \in \left[0, T_f\right]$, each point $X \in \Gamma(t)$ with a
well-defined normal moves in the direction normal to the front with a given
speed $F(X, \Gamma,t)$. Notice that $F$ may depend on the position, on the
time and on local properties of the front itself (certainly the normal
direction, not always defined, and possibly the local curvature or other
properties).

The problem is to approximate $\Gamma: \left[0, T_f\right] \rightarrow
\ensR^N$, given $\Gamma_0$ and $F$.

\subsubsection{Strategy}

The main idea is to evolve a function $\varphi: \ensR^N \times \left[0,
  T_f\right] \rightarrow \ensR$ such that
\begin{equation}
  \Forall t \in \left[0, T_f\right] \qquad \Gamma(t) = \left\{ x \in
    \ensR^N \Big{/} \varphi(x, t) = 0 \right\}.
  \label{eq:level_set_def}
\end{equation}

$\varphi$ is called the level set function. At any time, the zero level set of
$\varphi$ is the front itself. {\it A priori}, $\varphi$ could be any function
satisfying equation~(\ref{eq:level_set_def}). However, some assumptions (e.g.,
smoothness) and practical issues (e.g., initialization of $\varphi$) make it
convenient to define $\varphi$ as the signed distance to the front.

Then, if $d$ is the Euclidean distance on $\ensR^N$, we define, for any given
curve $\Upsilon$, the distance $d_{\Upsilon}$ to $\Upsilon$:
\begin{equation}
  \Forall x \in \ensR^N \qquad d_{\Upsilon}(x) = \min \left\{ d(x, P) \Big{/} P
    \in \Upsilon \right\}.
\end{equation}

Hence the signed distance $\varphi$ for all $x\in\ensR^N$ and $t\in\left[0,
  T_f\right]$:
\begin{equation}
  \varphi(x, t) = \left\{
    \begin{array}{ll}
      d_{\Gamma(t)}(x) & \textrm{if $x$ lies outside the front
        $\Gamma(t)$}\\ 
      - d_{\Gamma(t)}(x) & \textrm{if $x$ lies inside the front
        $\Gamma(t)$}\\ 
    \end{array}
  \right..
  \label{eq:phi}
\end{equation}

It can be shown that $\varphi$ obeys the equation
\begin{equation}
  \Forall x \in \ensR^N\quad \Forall t \in [0, T_f] \qquad
  \varphi_t(x, t) + F(x,\varphi(\cdot, t), t) {\| \nabla_x \varphi(x,
    t) \|}_2 = 0,
  \label{eq:level_set}
\end{equation}
where the velocity $F$ is now defined everywhere in $\ensR^N$ and depends on
the front through its dependence upon $\varphi$. Details may be found in
\citet{sethian99level}.

Recall that $\varphi(\cdot, 0)$ is known as well as $\Gamma_0$; $\varphi(0)$
is the signed distance to $\Gamma_0$:
\begin{equation}
  \Forall x\in\ensR^N \quad \varphi(x, 0) = \left\{
    \begin{array}{ll}
      d_{\Gamma_0}(x) & \textrm{if $x$ lies outside the front
        $\Gamma_0$}\\ 
      - d_{\Gamma_0}(x) & \textrm{if $x$ lies inside the front
        $\Gamma_0$}\\ 
    \end{array}
  \right..
  \label{eq:initial_cond}
\end{equation}

Equations~(\ref{eq:level_set}) and~(\ref{eq:initial_cond}) define the
initial-value problem that is to be solved. Zero level sets of $\varphi$ yield
the front points.

This nonstationary problem involves the Hamilton-Jacobi
equation~(\ref{eq:level_set}). There may be multiple solutions to this
equation. P.-L.  Lions and M. G. Crandall defined the so-called ``viscosity
solution'' of Hamilton-Jacobi equations \citep{lions82generalized,
  crandall83viscosity}, which turns out to be the unique physical solution for
which we search. Under given assumptions (mainly on the speed function $F$),
existence and uniqueness of the viscosity solution of the
problem~(\ref{eq:level_set})--(\ref{eq:initial_cond}) can be proved.

\subsection{Advantages and disadvantages of level set methods}

Several methods may be relevant to simulate the propagation of firefronts.
One may want to use marker techniques, in which the front is discretized by a
set of points. At each time step, each point is advanced according to the
speed function. This Lagrangian methodology leads to low-cost computations,
but requires care in the handling of topological changes.

Volume-of-fluid methods represent the front by the amount of each grid-cell
that is inside the front. In each cell, the front is approximated by a
straight line (horizontal or vertical, in most methods). Such methods can deal
with topological changes, but the front representation can be inaccurate. In
wildland firespread, the direction normal to the front is crucial because of
the wind-direction-dependent speed function (see
Section~\ref{sec:fire_model}).

Level set methods automatically deal with topological changes that occur in
wildland firespread, such as fronts merging and front convergence (in
connection with unburnt ``islands''). The level set description enables a fair
estimate of the normal to the front, making it well suited to the fire
propagation problem.
\\

However, level set methods have disavantages. First, they embed the front in a
higher-dimensional space. Helpfully, the narrow band level set method
\citep{adalsteinsson95fast} is an efficient algorithm which almost decreases
the problem dimension by one. Moreover, when it can be used, the fast marching
method \citep{sethian96fast} provides a highly efficient algorithm.

The main reservation may be the lack of proof of convergence of numerical
schemes for certain problems. For a given class of speed functions, the
problem~(\ref{eq:level_set})--(\ref{eq:initial_cond}) may routinely be solved
numerically \citep{crandall84two}. However, no proof of convergence in mesh
parameter or time step is yet available for some situations.


\subsection{Quick review of numerical approximations}

Numerical approximation to solutions of Hamilton-Jacobi equations is closely
related to numerical approximation to hyperbolic conservation
laws\footnote{Notice that, from equation~(\ref{eq:level_set}), $\varphi_x$
  satisfies a hyperbolic conservation law in the one-dimensional case.}. The
point is to introduce a numerical Hamiltonian to approximate the Hamiltonian
$H = F \cdot{\| \nabla_x \varphi \|}_2$.

Crandall and Lions have proven that, for given Hamiltonians and initial
conditions, a consistent, monotonic and locally Lipschitzian numerical
Hamiltonian yields a solution that converges to the viscosity solution.
Formal results may be found in \citet{crandall84two} and
\citet{souganidis85approximation}.

In one dimension, $\varphi_t + H(\nabla_x \varphi)=0$ may lead to the
following approximation:
\begin{equation}
  \varphi^{n+1}_{j} = \varphi^n_{j} - \Delta t \hspace{1mm} g
  \left(\frac{\varphi_{j+1} - \varphi_{j}}{\Delta x}, \frac{\varphi_{j} -
      \varphi_{j-1}}{\Delta x} \right).
  \label{eq:approx}
\end{equation}

For instance, if the Hamiltonian is not convex, the Lax-Friedrichs scheme may
be used; then, the numerical Hamiltonian is
\begin{equation}
  \Forall a, b \in \ensR \qquad g(a, b) = H\left(\frac{a+b}{2}\right)
  - \vartheta \frac{b-a}{2},
  \label{eq:Lax-Friedrichs}
\end{equation}
where the monotonicity is satisfied on $[-R, R]$ if $\displaystyle\vartheta =
\max_{-R \leq a \leq R} |H'(a)|$.

Several schemes have been developed, from simple and efficient schemes as that
of Engquist-Osher to high-order essentially nonoscillatory schemes
\citep{osher91high}.

\subsection{Overview of complexity issues}

Let the mesh (in $\ensR^N$) be orthogonal with $M$ points along each
direction. Assume that the front is described by $\mathcal{O}(M^{N-1})$
points. The narrow band level set method makes it sufficient to update the
level set function only in a narrow band (of width $k$) around the front. For
each time step, the complexity of the algorithm is therefore
$\mathcal{O}(kM^{N-1})$.

For an explicit temporal discretization the number of iterations is related to
the Courant-Friedrichs-Lewy condition.  Along $x$, the Courant number must be
less than $1$:
\begin{equation}
  \frac{\max \left|H'\right|\Delta  t}{\Delta x} \leq 1.
  \label{eq:CFL}
\end{equation}

Usually, controlling the accuracy of approximation is subordinate to space
discretization, which means that the time step is adjusted so that the Courant
number is taken close to $1$.
\\

Calculations may sometimes be sped up by reformulating the level set problem
as a stationary problem. This leads to the so-called fast marching method
\citep{sethian96fast}. Nevertheless, restrictions on the Hamiltonian prevent
the use of this technique for some applications.  The work of Sethian and
Vladimirsky has overcome some limitations \citep{sethian01ordered}, but
restrictive conditions still remain (e.g., convexity of the Hamiltonian).

\section{Code}
\label{sec:code}



\subsection{Introduction to the Multivac level set package}

Multivac is a level set package freely available (under the GNU GPL license)
at \url{http://vivienmallet.net/fronts/}. It is designed to be both efficient
and extensible, so that it may be used for a large range of applications.  To
achieve these goals, Multivac is built as a fully object-oriented library in
C++.

Multivac was designed independently of the firespread application described
herein, but easily enabled firespread simulations, and is presently
distributed with firespread-motivated functions. It has also been used in
modeling the growth of Si-based nanofilms \citep{phan03modeling} and image
segmentation.

The latest stable version available at the time of submission is
Multivac~1.10.

\subsection{Structure}

The modularity of Multivac comes from its object-oriented framework, in which
the main components of a simulation have been split into an equal number of
objects. A simulation is defined by the following objects:
\begin{itemize}
\item the \textit{mesh};
\item the \textit{level set function};
\item the \textit{velocity}, which provides the propagation rate of the front
  according to its position, its normal, its curvature, and the time;
\item the \textit{initial front};
\item the \textit{initializer}, which manages first initializations and
  initializations required by level set methods (e.g., the narrow band
  reconstruction);
\item the \textit{numerical scheme}, which advances the front in time;
\item the \textit{output management}.
\end{itemize}

For each item, a set of classes\footnote{A class is a user-defined type, in
  the manner of structures in C. Classes encapsulate data (called attributes)
  and functions (called methods).} with a common interface is available.  For
instance, several speed (i.e., propagation rate) functions are available
through several classes, e.g.  \code{CConstantSpeed} or \code{CFireModel}. All
speed functions have the same interface, which allows users to define their
own speed function on the same basis. The user principally provides speed
rates as a function of the position, the time, the normal to the front and the
curvature (these values are computed by Multivac itself).

\subsection{Calling sequence}

The whole is managed by an object of the class \code{CSimulator}. This object
simply calls the \textit{initializer} to perform the first initializations.
Then it manages the loop in time (or iterations, in the case of the fast
marching method) into which the \textit{numerical scheme} is called to advance
the front. The \textit{initializer} is called again to reinitialize the signed
distance function for the new step, and the object dedicated to
post-processing requirements is called to save any needed data.

In each step, objects communicate with one another through methods (i.e.,
functions) of their interface. For example, the \textit{velocity} object
provides speed rates to the \textit{numerical scheme}.

\subsection{Overview of available classes}

Multivac package (version 1.10) includes several classes which are listed in
Table~\ref{tab:classes}.

\begin{table}[htbp]
  \begin{tabular}{|c|c|}
    \hline
    Category & Available classes \\
    \hline
    \hline
    Mesh & Orthogonal mesh \\
    \hline
    Level set function & Defined on an orthogonal mesh \\
    \hline
    Velocity & Constant speed \\
    & Piecewise constant speed \\
    & Fire model \\
    & Simplified fire model \\
    & Image intensity \\
    & Image gradient \\
    \hline
    Initial front & Circle \\
    & Two or three circles \\
    & One or two circles with an island inside \\
    & Front defined by any set of points \\
    \hline
    Initializer & Basic initialization (no velocity extension) \\
    & Extends the velocity with the closest \\
    & neighbor on the front \\
    \hline
    Numerical scheme & Engquist-Osher, first order \\
    (narrow band) & Lax-Friedrichs, first order \\
    & Engquist-Osher, ENO, second order \\
    & Chan-Vese algorithm \citep{chan01active} \\
    \hline
    Numerical scheme & Engquist-Osher, first order \\
    (fast marching) & \\
    \hline
  \end{tabular}
  \caption{Basic classes available in Multivac 1.10.}
  \label{tab:classes}
\end{table}

\subsection{Other strengths, limitations and future work}

Multivac takes advantage of C++ exceptions to track errors, and several
debugging levels are defined, from a safe mode, in which all is checked, to a
fast mode, in which performance is the primary concern.

There are currently two main limitations. First, Multivac deals only with
uniform orthogonal meshes. However, extensions of level set methods to
unstructured meshes exist (e.g., \citet{barth98numerical}) and they could be
implemented within the Multivac framework. Adaptively refined meshes are also
accommodated with additional mathematical complexity, though the
implementation effort would be substantial. Second, Multivac deals only with
two-dimensional problems.

Work is planned to allow inverse modeling (parameter estimation based on data
assimilation) within the framework of Multivac. The main idea is to replace
the class \code{CSimulator} with a class dedicated to inverse modeling.
Preliminary results show the framework extendibility, but this capability is
not yet available in distributed versions. Future versions should include this
feature, based on an innovative method for integrating sensitivities along
with the front itself.

\section{Complexity and Convergence Studies}
\label{sec:complexity}



\subsection{Convergence studies}

In this section, we report convergence studies that are necessary to validate
the code. As in \citet{adalsteinsson98fast}, tests are carried out for a
circle that expands in time with a unitary velocity. Details of the simulation
are summarized in Table~\ref{tab:test_case}.

\begin{table}[htbp]
  \begin{tabular}{|c|c|c|}
    \hline
    Data & Value & Comment\\
    \hline
    Domain & $\Omega = [0, 3] \times [0, 3]$ & \\
    Initial front & Circle & \\
    Circle center & $\left(x_c, y_c\right) = \left(1.5, 1.5\right)$ & Domain
    center\\
    Initial circle radius & $r_{initial} = 0.5$ & \\
    Final circle radius & $r_{final} = 0.9$ & \\
    Velocity & $F = 1.0$ & Constant\\
    Duration & $T_f = 0.4$ & \\
    Time step & $\Delta t = 10^{-4}$ & \\
    \hline
  \end{tabular}
  \caption{Simulation test-case.}
  \label{tab:test_case}
\end{table}

We introduce three norms. The first is
\begin{equation}
  \label{eq:es}
  e_{spatial}^1 = |r_{simulated} - r_{final}|,
\end{equation}
where $r_{simulated}$ is the simulated radius, estimated as follows:
\begin{equation}
  r_{simulated} = \frac{1}{\mathrm{card}(\Gamma_d)} \sum_{(x, y)\in\Gamma_d}
  \textrm{d} \left( (x, y), \left(x_c, y_c\right) \right),
\end{equation}
where $\Gamma_d$ is the discretized front as returned by the simulation (at
time $T_f$) and $\textrm{d}$ is the Euclidian distance.

Additionally, if $T_{true}(x, y)$ is the time at which the front is supposed
to reach the point $(x, y)$:
\begin{equation}
  \label{eq:et}
  e_{time}^2 = \sqrt{ \frac{1}{\mathrm{card}(\Gamma_d)} \sum_{(x,
      y)\in\Gamma_d} \left( T_f - T_{true}(x, y) \right)^2 }.
\end{equation}

The last norm is an infinity norm:
\begin{equation}
  \label{eq:ef}
  e_{time}^\infty = \max_{(x, y)\in\Gamma_d} \left| T_f - T_{true}(x,
    y) \right|.
\end{equation}

Table~\ref{tab:eo} shows results for the first-order Engquist-Osher scheme
with the narrow band method. The width of the band is 12 cells and the front
lies within a band whose width is 6 cells.

\begin{table}[htbp]
  \begin{tabular}{|c|c|c|c|c|}
    \hline
    $\Delta x = \Delta y$ & $N_x = N_y$ & $e_{spatial}^1$
    $\left(\times 10^3\right)$ & $e_{time}^2$ $\left(\times
      10^3\right)$ & $e_{time}^\infty$ $\left(\times 10^3\right)$ \\
    \hline
    0.01 & 301 & 1.634 & 1.753 & 2.377 \\
    0.005 & 601 & 0.855 & 0.901 & 1.191 \\
    0.0025 & 1,201 & 0.460 & 0.474 & 0.600 \\
    0.00125 & 2,401 & 0.244 & 0.247 & 0.299 \\
    \hline
  \end{tabular}
  \caption{Errors versus spatial discretization.}
  \label{tab:eo}
\end{table}

The first-order Lax-Friedrichs scheme and the second ENO Engquist-Osher scheme
were also checked successfully. As for the second-order scheme, the
full-matrix method, that is, without the narrow-band restriction, was used
because the front reconstruction destroys the second-order accuracy.

\subsection{Complexity issues}

Multivac was compiled under Linux with GNU/g++ 3.3, and the reference
simulation (see Table~\ref{tab:test_case}) was launched on a Pentium 4 running
at 2.6 Ghz.  The width of the narrow band was 12 cells and the width of the
inner band, in which the front lies, was 6 cells. If $N_x=N_y=1001$ (one
million cells), the $4,000$ iterations were achieved in 14\ s.

The complexity of the narrow band level set method is close to
$\mathcal{O}(N)$, where $N=N_x=N_y$. Table~\ref{tab:timings} shows that linear
complexity of the method is not observed. Instead, the complexity seems to be
$\mathcal{O}(N^2)$. This is the complexity of the suboptimal algorithm
currently used to rebuild the front. Moreover, the number of front
reconstructions increases with the mesh refinement since the width of the
narrow band does not change.

\begin{table}[htbp]
  \begin{tabular}{|c|c|c|c|c|}
    \hline
    $\Delta x = \Delta y$ & $N_x = N_y$ & Timings (s) \\
    \hline
    0.03 & 101 & 0.4 \\
    0.015 & 201 & 0.9 \\
    0.01 & 301 & 1.6 \\
    0.0075 & 401 & 2.6 \\
    0.006 & 501 & 4.0 \\
    0.005 & 601 & 5.6 \\
    0.004285714 & 701 & 7.4 \\
    0.00375 & 801 & 9.5 \\
    0.003333333 & 901 & 11.9 \\
    0.003 & 1001 & 14.1 \\
    \hline
  \end{tabular}
  \caption{Timings versus spatial discretization.}
  \label{tab:timings}
\end{table}

\section{Applying Level Set Methods to Firespread Applications}
\label{sec:applying}




\subsection{Method and numerical scheme}

The speed function~(\ref{eq:fire_model}) introduced in the level set
equation~(\ref{eq:level_set}) provides an Hamiltonian with nontrivial
dependencies.  Because of these dependencies (particularly the non-convexity
of the Hamiltonian), neither the fast marching method nor its extension to
anisotropic problems can be applied. The narrow-band level set method is more
relevant.

A highly accurate numerical scheme is not required for the investigations
reported here. The discrepancies between the numerical simulation and the
exact solution should be considered in the context of other approximations:
the model itself is simplistic; input parameters such as wind speed or fuel
density are typically not accurately estimated; the location of the initial
front introduces further uncertainties. A first-order scheme suffices for our
purposes.

Since the Hamiltonian involved is not convex with respect to spatial
derivatives of the level set function, the first-order Lax-Friedrichs scheme
(refer to equation~(\ref{eq:Lax-Friedrichs})) is well suited. To minimize
introduction of diffusivity, a local Lax-Friedrichs scheme may be used as
well.

As previously advocated, the timestep $\Delta t$ is chosen according to the
Courant-Friedrichs-Lewy condition~(\ref{eq:CFL}):
\begin{equation}  
  \Delta t = \frac{\alpha\Delta x}{\max \left|H'\right|},
\end{equation}
where $\alpha \leq 1$; $\alpha$ is not kept constant in the tests that we
undertake. Nevertheless, the Courant-Friedrichs-Lewy condition is estimated at
every iteration with an ({\it a priori}) approximation to $\max
\left|H'\right|$ along $x$ and $y$, which leads to:
\begin{equation}
  \Delta t \leq \frac{\Delta x}{a (m + 1) \sqrt{U}}.
\end{equation}

The main characteristics of the simulation, including model parameters (refer
to equation~(\ref{eq:fire_model})), are gathered in Table~\ref{tab:param}.

\begin{table}[htbp]
  \begin{center}
    \subtable{
      \begin{tabular}{|c|c|}
        \hline
        Parameter & Value\\
        \hline
        \hline
        $n$ & 1.5 \\
        $U$ & 100 \\
        $a$ & 0.5 \\
        $\varepsilon_0$ & 0.2 \\
        $\alpha$ & 0.5 \\
        \hline
        \hline
        $\Delta x$ & $3\cdot 10^{-3}$ \\
        $\Delta y$ & $3\cdot 10^{-3}$ \\
        $\Delta t$ & $10^{-4}$ \\
        $T_f$ & 0.1 \\
        \hline
      \end{tabular} }%
    \subtable{
      \begin{tabular}{|c|c|}
        \hline
        Parameter & Value\\
        \hline
        \hline
        Domain & $\Omega = [0, 3] \times [0, 3]$ \\
        Initial front & Circle \\
        Circle center & $\left(1.5, 1.0\right)$ \\
        Initial circle radius & $r_{initial} = 0.5$ \\
        Velocity & $F = 1.0$ \\
        Duration & $T_f = 0.1$ \\
        Time step & $\Delta t = 5\cdot 10^{-5}$ \\
        Spatial discretization & $N_x = N_y = 1001$ \\
        \hline
      \end{tabular} }%
    
    \caption{Parameters and their default values.}
    \label{tab:param}
  \end{center}
\end{table}

\subsection{Results}

\begin{figure}
  \includegraphics[width=0.8\textwidth]{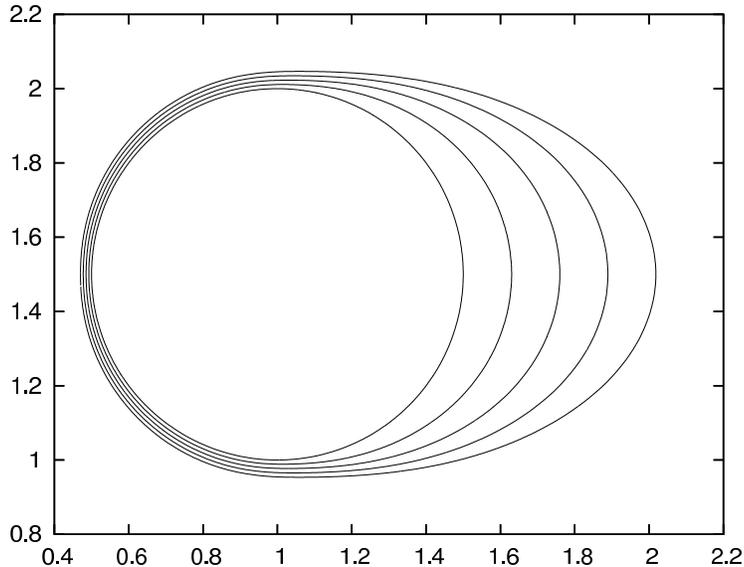}
  \caption{Basic simulation described by Table~\ref{tab:param}.}
  \label{fig:basic_case}
\end{figure}

The simulation described by Table~\ref{tab:param} is shown in
Figure~\ref{fig:basic_case}. The figure shows snapshots of the front,
initially circular, at subsequent times, under a constant-magnitude wind
blowing from left to right. Since thoroughly burnt areas cannot be burnt again
(on the time scale of the simulation), the area enclosed by the front
increases with time. The rear, the flanks and the head of the front are
clearly identifiable.

\begin{figure}
  \includegraphics[width=0.8\textwidth]{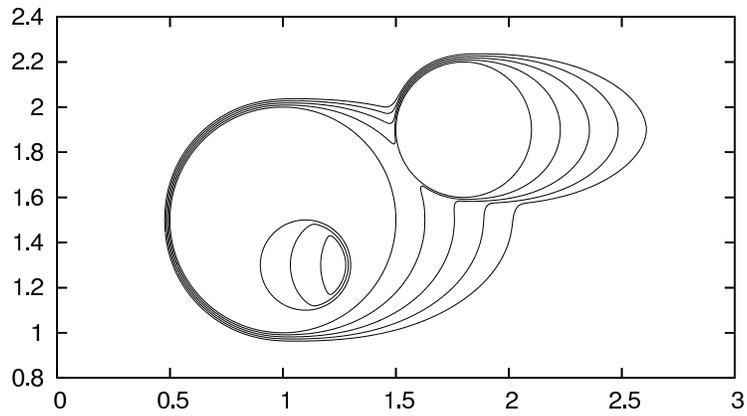}
  \caption{Two main fronts merge, and an island -- the unburnt area within the
    biggest front -- is burnt.}
  \label{fig:merge_and_island}
\end{figure}

The reference simulation is slightly modified to show the ability to deal with
multiple fronts -- Figure~\ref{fig:merge_and_island}. It demonstrates the
capability to deal with the merging of fronts (two main fronts), and to deal
with the so-called islands, i.e. an unburnt area surrounded by a burnt area.

\begin{figure}
  \includegraphics[width=0.8\textwidth]{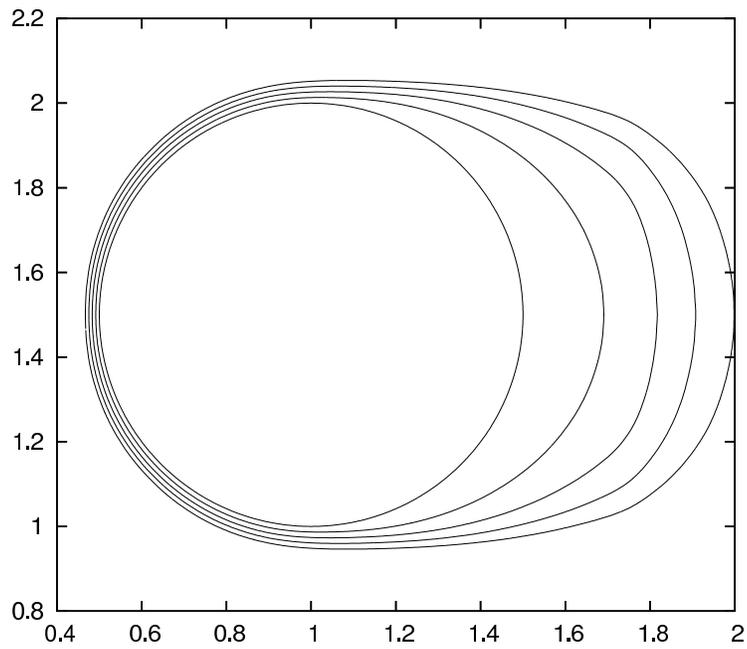}
  \caption{The front slows down at the head for $a = 0.25$ if $x>1.8$. The
    final time is changed to $T_f = 1.5$.}
  \label{fig:fuel_slow}
\end{figure}

\begin{figure}
  \includegraphics[width=0.8\textwidth]{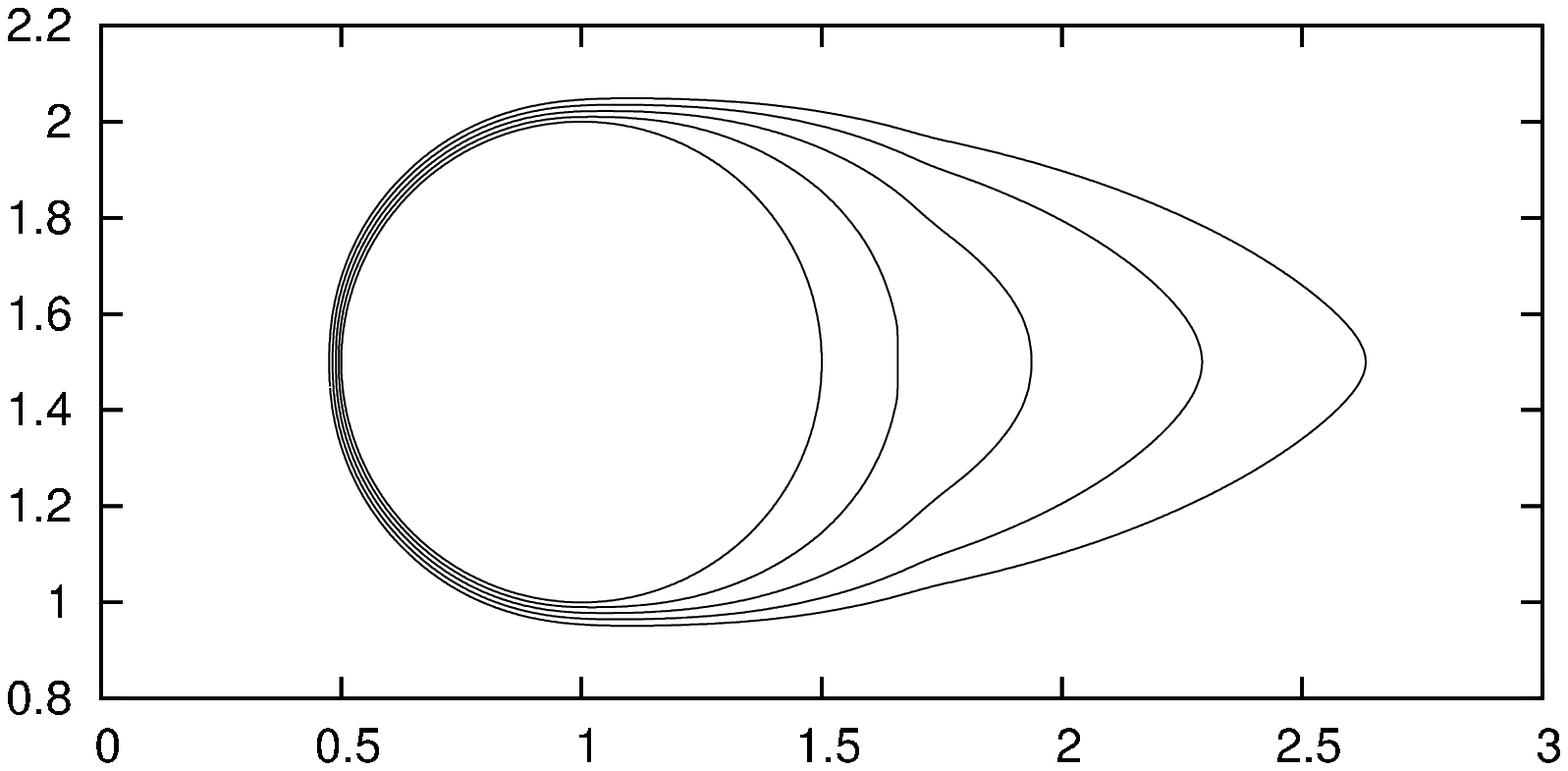}
  \caption{The front advances faster at the head for $a = 1.0$ if $x>1.8$.}
  \label{fig:fuel_fast}
\end{figure}

In Figures~\ref{fig:fuel_slow} and~\ref{fig:fuel_fast}, we use the same
parameters as in Table~\ref{tab:param} but $\Delta t = 2.5\cdot 10^{-5}$, and
$a$ depends on $x$, $a$ being equal to $0.5$ if $x<1.7$, and $a=0.25$
(Figure~\ref{fig:fuel_slow}) or $a=1.0$ (Figure~\ref{fig:fuel_fast}) if
$x>1.8$, and $a$ being linearly interpolated for intermediate values of $x$.
Since $a$ takes into account the available fuel loading, these two simulations
roughly show the influence of the inhomogeneous available fuel loading, should
it increase (Figure~\ref{fig:fuel_slow}) or decrease
(Figure~\ref{fig:fuel_fast}). The inherent decrease of the radius of curvature
at the head for a constant-direction wind suggests that some vacillation of
wind direction contributes when the head broadens under otherwise uniform
conditions.

\begin{figure}
  \includegraphics[width=0.8\textwidth]{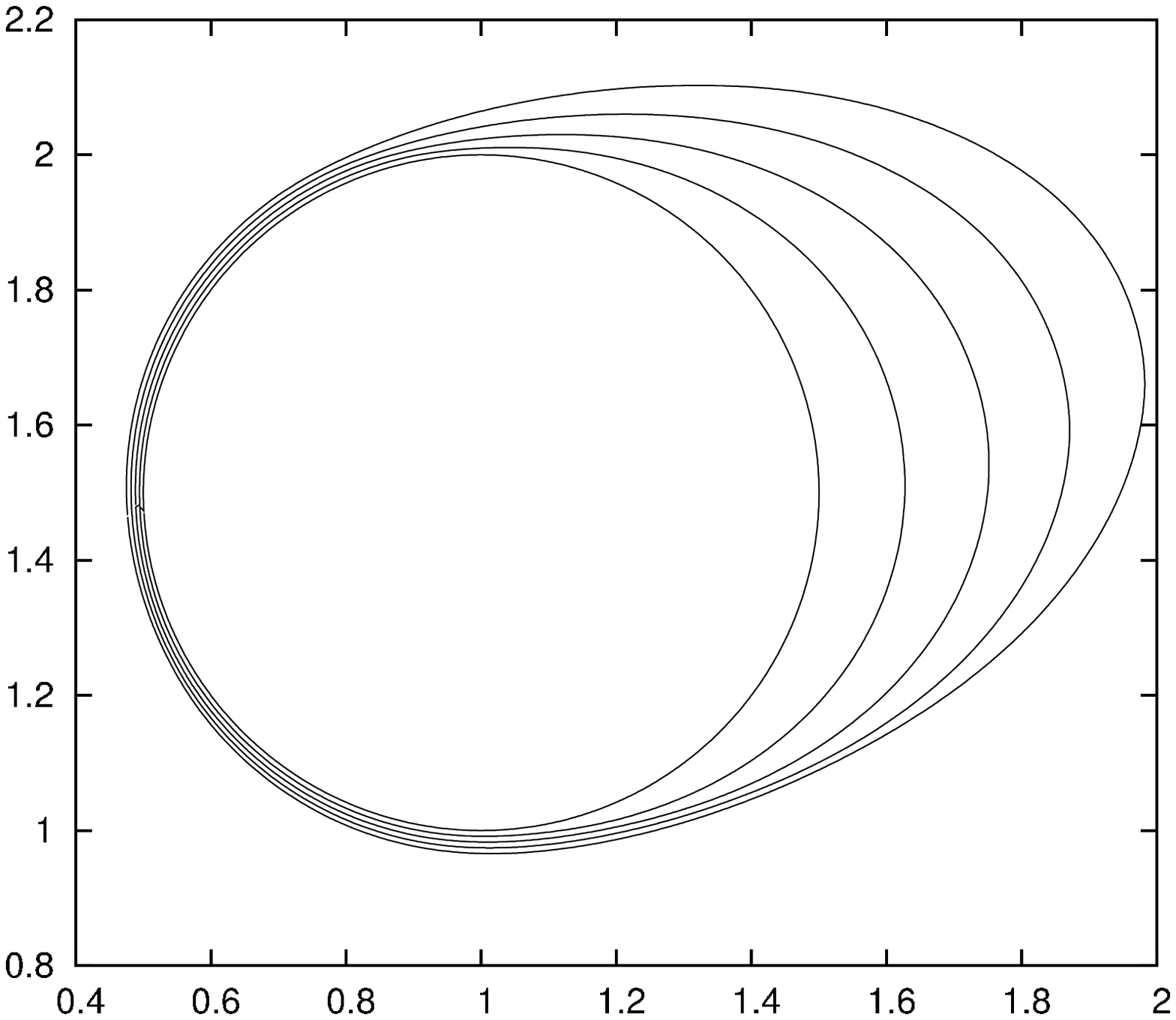}
  \caption{Same as the reference simulation, but with a changing wind
    direction.}
  \label{fig:rotating}
\end{figure}

Figure~\ref{fig:rotating} shows the impact of a rotating wind direction. If
north is toward the top of the figure, then the wind is oriented first
west-to-east and tends later to south-to-north.

\begin{figure}
  \includegraphics[width=0.8\textwidth]{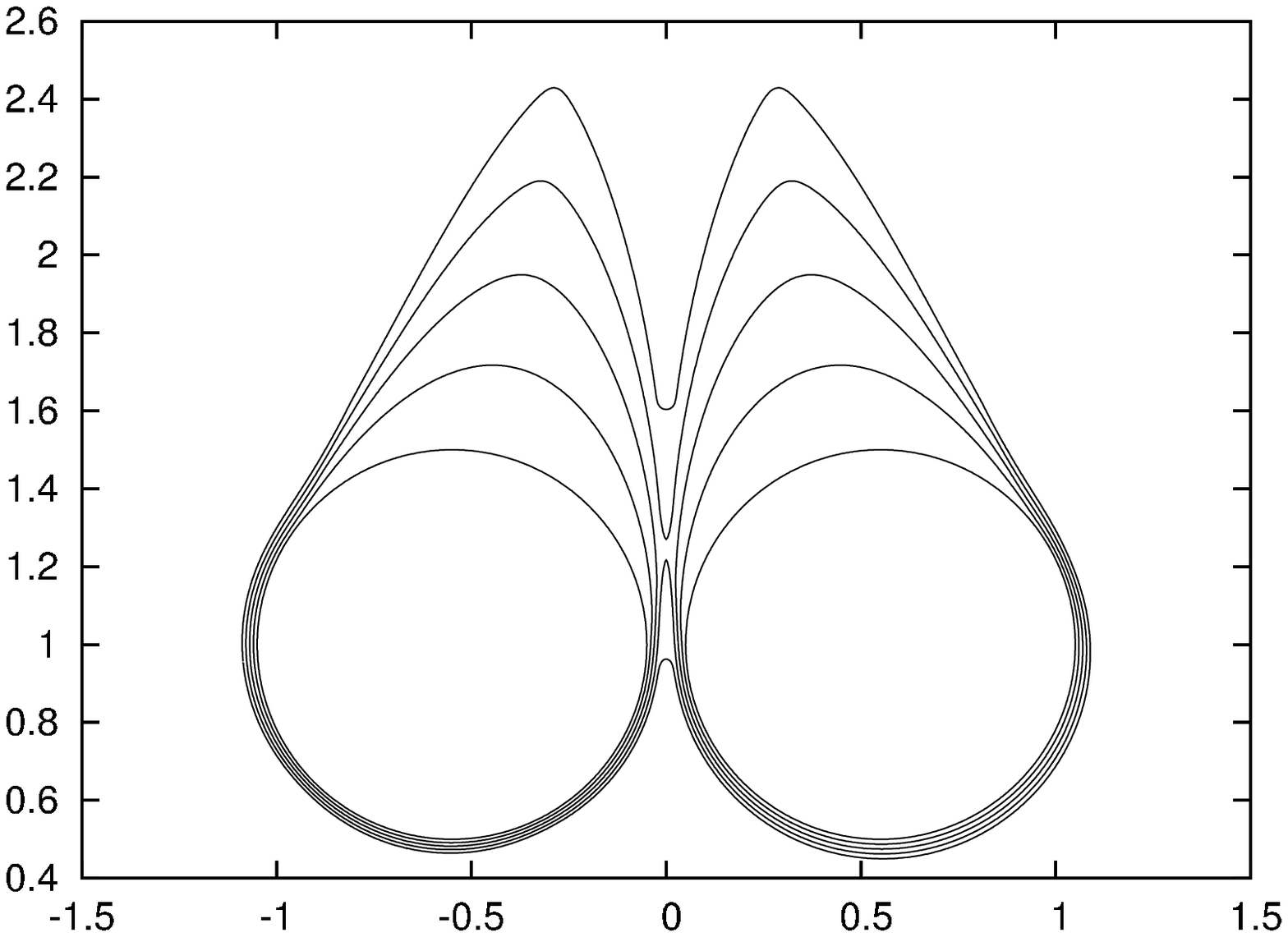}
  \caption{Evolution of the merged front from initially two mirror-image
    fronts, one to each side of the stagnation line for a converging
    $x$-component wind, but both to one side of the stagnation line for a
    diverging $y$-component wind.}
  \label{fig:counterflow_up}
\end{figure}

\begin{figure}
  \includegraphics[width=0.8\textwidth]{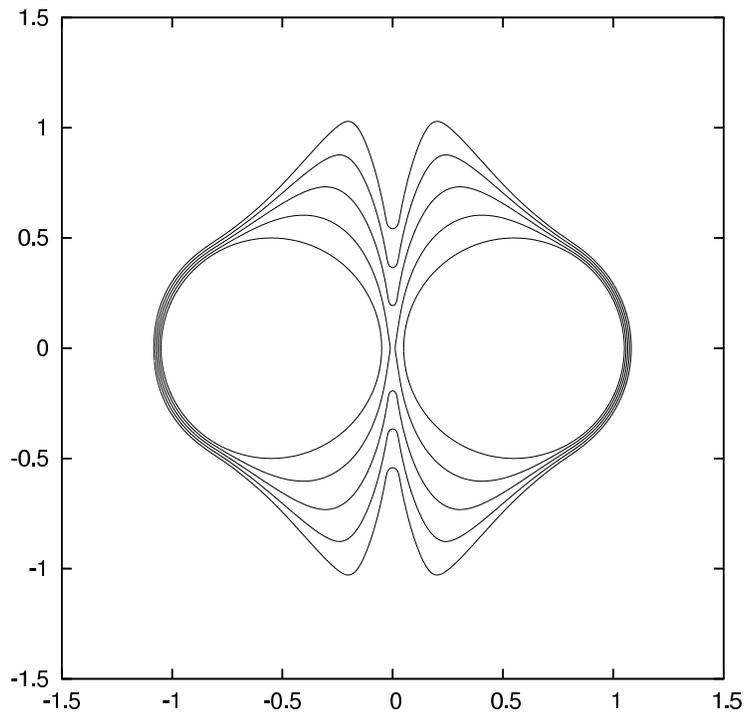}
  \caption{Evolution of the merged front from initially two mirror-image
    fronts, here symmetrically sited relative to a simple counterflow wind.}
  \label{fig:counterflow_middle}
\end{figure}

The next two Figures~\ref{fig:counterflow_up} and~\ref{fig:counterflow_middle}
show the behavior of two fronts subject to a simple-counterflow wind, i.e., a
wind defined as:
\begin{equation}
  \overrightarrow{U}(x, y) = \left(
    \begin{array}{c}
      -u x\\ 
      u y 
    \end{array}
  \right)
\end{equation}
where $u$ is set to $100$. A counterflow exemplifies wind conditions well
suited for setting a backfire, to preburn the vegetation in the path of a
wind-aided fire.

\begin{figure}
  \includegraphics[width=0.8\textwidth]{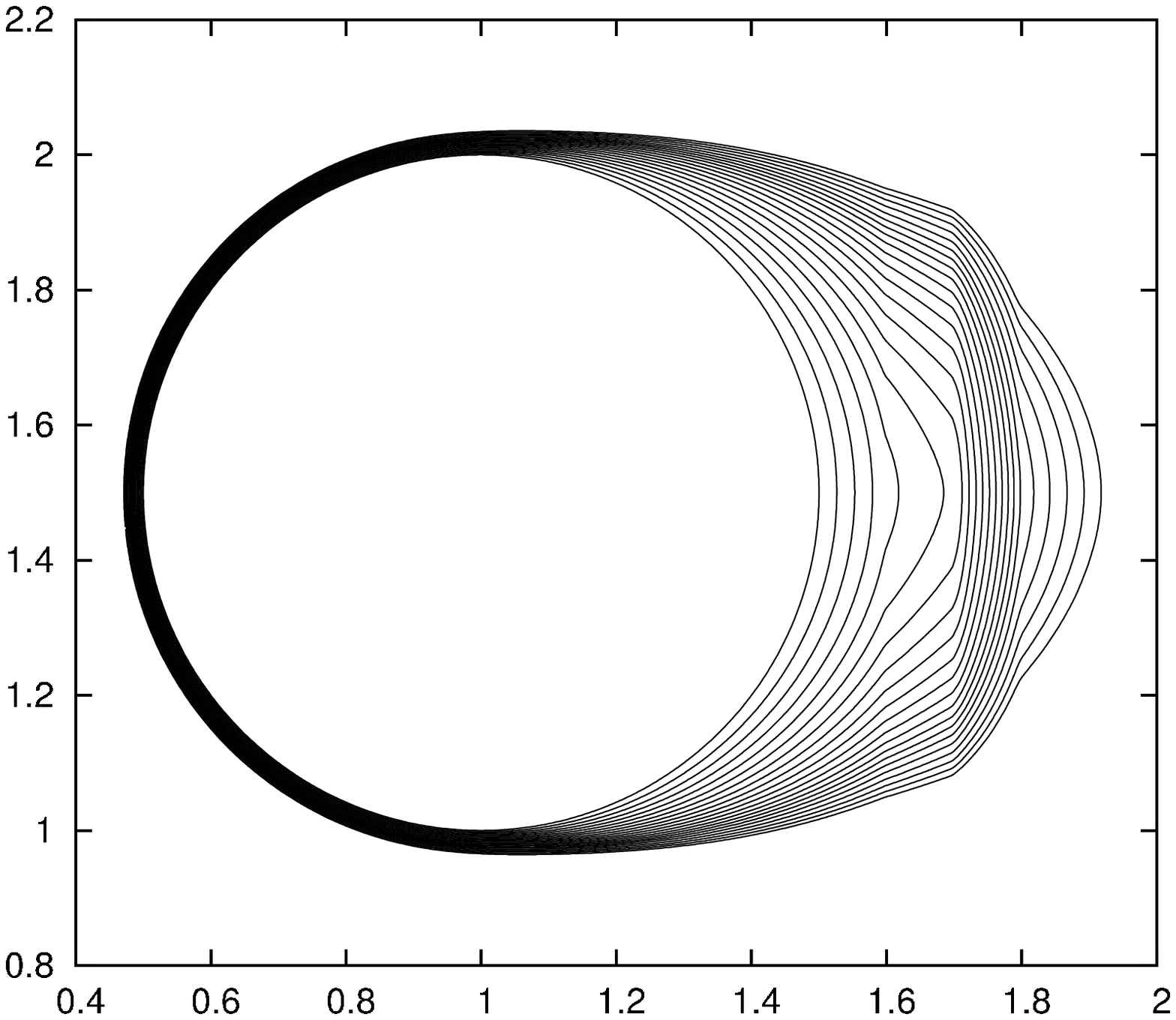}
  \caption{Taking into account topography: the front propagates over an
    idealized hill.}
  \label{fig:topography}
\end{figure}

The last Figure~\ref{fig:topography} shows a front that propagates over an
idealized hill. Where the slope is positive (between $x=1.6$ and $x=1.7$), the
firefront typically advances faster.  Downhill the front typically slows down
\citep[][pp. 94--97]{luke78bushfires}.  The speed function is therefore
modified to take into account the slope $s$:
\begin{equation}
  \label{eq:topography}
  F_{\mathrm{topography}} = F \times e^{2s},
\end{equation}
where $s$ is in radians.

\section{Conclusion and Future Prospects}
\label{sec:conclusion}

A semi-empirical, equilibrium-type firespread rate has been used to model a
wind-aided firefront propagation across wildland surface vegetation.  In this
formulation, the rate depends primarily on the wind speed, and the angle
between the wind direction and the normal to the firefront (idealized as a
one-dimensional interface). In scenarios arising in practice, the front may
consist of several closed curves (possibly nested) that can merge as they
propagate.

Level set methods appear capable of treating the model formulated to simulate
wildland fire evolution.  They treat readily the topological changes that may
occur to the firefront, and they are known to converge to the physical
solution of front tracking problem.

They were applied via the Multivac package.  This open-source library is
designed to handle a wide range of applications without loss of computing
performance. It includes several algorithms and numerical schemes, primarily
for the narrow-band level set method, which is more computationally efficient
than the full level set method.

A possible direction for future work is to focus on parameter estimation
within the context of the simple model illustrated herein. A cost is
introduced to measure the distance between the simulated front and ground,
aerial, and/or satellite observations. The discrepancy between the simulated
and observed positions of the front may be based either on the front arrival
times (at monitored locations), or on distances between the simulated front
and the monitored locations (at arrival times). For gradient-based
optimization methods, the main challenge is to compute the derivative of the
cost function with respect to the parameters. An adjoint code being difficult
to construct, alternative methods should be sought.

This work could help guide fire-control tactics. The objective function would
then penalize front advance into societal assets, and penalize the cost of the
firefighting activity. The parameters would be the model variables modifiable
by firefighting countermeasures.  The links between this optimization problem
and shape optimization should be investigated.








\section*{Acknowledgments}

The support of the National Science Foundation under grant CCF-03-52334 and
the U.S.  Department of Agriculture Forest Service under grant SFES
03-CA-11272169-33, administered by the Riverside Forest Fire Laboratory, a
research facility of the Pacific Southwest Research Station, is gratefully
acknowledged.  The authors are particularly indebted to Dr.~Francis M.~Fujioka
of Riverside for enhancing the relevance of our research through his technical
advice, and for his support for the training of summer students in the
computational technology of firespread and fire imaging.

\bibliographystyle{apalike} \bibliography{references}

\end{document}